\documentclass[11pt,a4paper,final]{article}
\usepackage[utf8]{inputenc}
\usepackage[english]{babel}
\usepackage{amsmath}
\usepackage{amsfonts}
\usepackage{amssymb}
\usepackage{graphicx}
\begin{document}

\title{On a higher order multi-term time-fractional partial differential equation involving Caputo-Fabrizio derivative}

\medskip
\author{E.T.Karimov, S.Pirnafasov}
\maketitle

\begin{abstract}
In the present work we discuss higher order multi-term partial differential equation (PDE) with the Caputo-Fabrizio fractional derivative in time. We investigate a boundary value problem for fractional heat equation involving higher order Caputo-Fabrizio derivatives in time-variable. Using method of separation of variables and integration by parts, we reduce fractional order PDE to the integer order. We represent explicit solution of formulated problem in particular case by Fourier series.   
\end{abstract}



\section{Introduction}
Consideration of new fractional derivative with non-singular kernel was initiated by M.Caputo and M.Fabrizio in their work \cite{CF1}. Motivation came from application. Precisely,  new fractional derivatives can better describe material heterogeneities and structures with different scales. Special role of spatial fractional derivative in the study of the macroscopic behaviors of some materials, related with nonlocal interactions, which are prevalent in determining the properties of the material was also highlighted. In their next work \cite{CF2}, authors represented some applications of the introduced fractional derivative. J.J.Nieto and J.Losada studied some properties of this freactional derivative naming it as Caputo-Fabrizio (CF) derivative \cite{NL}. Namely, they introduced fractional integral associated with the Caputo-Fabrizio derivative, applying it to the solution of linear and nonlinear differential equations involving CF derivative.  

Later, many authors showed interest to this CF derivative and as a result, several applications of CF derivative were discovered. For instance, in groundwater modeling \cite{AA1},\cite{AB1}, in electrical circuits \cite{AA2}, in controlling the wave movement \cite{AA3}, in nonlinear Fisher's reaction-diffusion equation \cite{A1}, in modeling of a mass-spring-damper system \cite{GA1} and etc. We note also other works related with CF derivative \cite{GA2}, \cite{GY}, \cite{ASh}, \cite{AN}, \cite{GA3}, \cite{H}.

Different methods were applied for solving differential equations involving CF derivative. Namely, Laplace transform, reduction to the integral equation and reduction to the integer order differential equation. Last two methods were used in the works \cite{K1}, \cite{K2}. 

In this paper, we aim to show algorithm how to reduce initial value problem (IVP) for multi-term fractional DE with CF derivative to the IVP for integer order DE and using this result to prove a unique solvability of a boundary value problem (BVP) for partial differential equation (PDE) involving CF derivative on time-variable. First, we give preliminary information on CF derivative and then we formulate our main problem. Representing formal solution of the formulated problem by infinite series, in particular case, we prove uniform convergence of that infinite series.

\section{Preliminaries}
\begin{equation}
{}_{CF}D_{at}^{\alpha}g(t)=\frac{1}{1-\alpha}\int\limits_a^t g'(s)e^{-\frac{\alpha}{1-\alpha}(t-s)}ds
\label{eqn2.1}
\end{equation}
is a fractional derivative of order $\alpha\, (0<\alpha<1)$ in Caputo-Fabrizio sense \cite{CF2}. We note that operator
(2) is well defined by the set
\[
W^{\alpha,1}=\left\{g(t)\in L^1(a,\infty);\,(f(t)-f_a(s))e^{-\frac{\alpha}{1-\alpha}(t-s)}\in L^1(a,t)\times L^1(a,\infty)\right\},
\]
whose norm is given by for $\alpha\neq 1$ by
\[
||g(t)||_{W^{\alpha,1}}=\int\limits_a^\infty |g(t)dt+\frac{\alpha}{1-\alpha}\int\limits_a^\infty\int\limits_{-\infty}^t
|g_s(s)|e^{-\frac{\alpha}{1-\alpha}(t-s)}ds dt,
\]
where $g_a(t)=g(t),\,t\geq a,\,\,\,g_a(t)=0,\,\,\,-\infty<t<a$ \cite{CF2}. Moreover, the following equality 
\[
{}_{CF}D_{at}^{\alpha+n}g(t)={}_{CF}D_{at}^\alpha\left({}_{CF}D_{at}^n g(t)\right)
\]
is true \cite{GA1}.
\section{Formulation of a problem and formal solution}
Consider the following time-fractional partial differential equation
\begin{equation}
\sum\limits_{n=0}^{k} \lambda_n \,\cdot\,{}_{CF}D_{0t}^{\alpha+n}u(t,x)-u_{xx}(t,x)=f(t,x)
 \label{eqn3.1}
\end{equation}
in a domain $\Omega=\left\{(t,x):\, 0<t<q,\,0<x<1\right\}$. Here $\lambda_n$ are given real numbers, $f(t,x)$ is a given function, $k\in \mathbf{N}_0,\,\,q\in \mathbf{R}^+$.

\textbf{Problem.} To find a solution of Eq.(\ref{eqn3.1}) satisfying the following conditions:
\begin{itemize}
\item 
\begin{equation}
u(t,x)\in C^2(\Omega),\,\,\frac{\partial^n u(t,x)}{\partial t^n}\in W^{\alpha,1}(0,q);
\label{eqn3.2}
\end{equation}

\item 
\begin{equation}
u(t,0)=u(t,1)=0,\,\, \left.\frac{\partial^i u(t,x)}{\partial t^i}\right|_{t=0}=\tilde{C}_i,\,\,i=0,1,2,...,k,
\label{eqn3.3}
\end{equation}
where $\tilde{C}_i$ are any real numbers.
\end{itemize}
\subsection{Solution of higher order multi-term fractional ordinary differential equation}
We expand $u(t,x)$ to the Fourier series as
\begin{equation}
u(t,x)=\sum\limits_{m=1}^\infty T_m(t)\sin m\pi x, \label{eqn3.4}
\end{equation}
where $T_m(t)$ are Fourier coefficients of $u(t,x)$.

Substituting representation (\ref{eqn3.4}) into Eq.(\ref{eqn3.1}) and considering initial conditions (\ref{eqn3.3}), we get the following initial value problem with respect to time-variable:
\begin{equation}
\left\{
\begin{array}{l}
\sum\limits_{n=0}^{k}\lambda_n\,\cdot\, {}_{CF}D_{0t}^{\alpha+n}T_m(t)+(m\pi)^2T_m(t)=f_m(t),\\
T_m^{(i)}(0)=C_i,\,\,i=0,1,2,...,k\\
\end{array}
\right.\label{eqn3.5}
\end{equation}
where $f_m(t)$ are Fourier coefficients of $f(t,x)$.

Based on definition (\ref{eqn2.1}) and initial conditions (\ref{eqn3.3}), after applying of the formula for integration by parts, we rewrite fractional derivatives as follows:
\[
\begin{array}{l}
{}_{CF}D_{0t}^{\alpha}T_m(t)=\frac{1}{1-\alpha}T_m(t)-\frac{\alpha}{(1-\alpha)^2}\int\limits_0^t T_m(s)e^{-\frac{\alpha}{1-\alpha}(t-s)}ds-\frac{T_m(0)}{1-\alpha}e^{-\frac{\alpha}{1-\alpha}t},\\
{}_{CF}D_{0t}^{\alpha+1}T_m(t)=\frac{1}{1-\alpha}T_m'(t)-\frac{\alpha}{(1-\alpha)^2}T_m(t)+
\frac{\alpha^2}{(1-\alpha)^3}\int\limits_0^t T_m(s)e^{-\frac{\alpha}{1-\alpha}(t-s)}ds+\\
\frac{\alpha T_m(0)}{(1-\alpha)^2}e^{-\frac{\alpha}{1-\alpha}t}-\frac{T_m'(0)}{1-\alpha}e^{-\frac{\alpha}{1-\alpha}t},\\
{}_{CF}D_{0t}^{\alpha+2}T_m(t)=\frac{1}{1-\alpha}T_m''(t)-\frac{\alpha}{(1-\alpha)^2}T_m'(t)+
\frac{\alpha^2}{(1-\alpha)^3}T_m(t)-\\\frac{\alpha^3}{(1-\alpha)^4}\int\limits_0^t T_m(s)e^{-\frac{\alpha}{1-\alpha}(t-s)}ds-
\frac{\alpha^2 T_m(0)}{(1-\alpha)^3}e^{-\frac{\alpha}{1-\alpha}t}+\frac{\alpha T_m'(0)}{(1-\alpha)^2}e^{-\frac{\alpha}{1-\alpha}t}-\frac{T_m''(0)}{1-\alpha}e^{-\frac{\alpha}{1-\alpha}t},\\
...\\
\end{array}
\]
Continuing this procedure, we can find for $n\geq 1$ the following formula
\begin{equation}
\begin{array}{l}
{}_{CF}D_{0t}^{\alpha+n}T_m(t)=\frac{1}{1-\alpha}\left\{\sum\limits_{i=0}^n \left(-\frac{\alpha}{1-\alpha}\right)^i
\left[ T_m^{(n-i)}(t)-T_m^{(n-i)}(0)e^{-\frac{\alpha}{1-\alpha}t}\right]+\right.\\
\left.\left(-\frac{\alpha}{1-\alpha}\right)^{n+1}\int\limits_0^t T_m(s)e^{-\frac{\alpha}{1-\alpha}(t-s)}ds\right\}.
\end{array}
\label{eqn3.6}
\end{equation}
We substitute (\ref{eqn3.6}) into (\ref{eqn3.5}) and deduce
$$
\begin{array}{l}
\sum\limits_{n=0}^k\frac{\lambda_n}{(1-\alpha)}\left\{\sum\limits_{i=0}^{n}\left(-\frac{\alpha}{1-\alpha}\right)^i \left[T_m^{(n-i)}(t)-T_m^{(n-i)}(0)e^{-\frac{\alpha}{1-\alpha}t}\right]+\right.\\
\left.\left(-\frac{\alpha}{1-\alpha}\right)^{n+1}\int\limits_0^t T_m(s)e^{-\frac{\alpha}{1-\alpha}(t-s)}ds\right\}+(m\pi)^2T_m(t)=f_m(t).
\end{array}
$$

We multiply this equality by $(1-\alpha)e^{\frac{\alpha}{1-\alpha}t}$:
\begin{equation}
\begin{array}{l}
\sum\limits_{n=0}^k \lambda_n\left\{\sum\limits_{i=0}^n\left(-\frac{\alpha}{1-\alpha}\right)^i \left[T_m^{(n-i)}(t)e^{\frac{\alpha}{1-\alpha}t}-T_m^{(n-i)}(0)\right]+\right.\\
\left.\left(-\frac{\alpha}{1-\alpha}\right)^{n+1}\int\limits_0^t T_m(s)e^{\frac{\alpha}{1-\alpha}s}ds\right\}+(m\pi)^2(1-\alpha)T_m(t)e^{\frac{\alpha}{1-\alpha}t}\\
=(1-\alpha)e^{\frac{\alpha}{1-\alpha}t}f_m(t).
\end{array}
\label{eqn3.7}
\end{equation}
Introducing new function as $\tilde{T}_m(t)=T_m(t)e^{\frac{\alpha}{1-\alpha}t}$, we rewrite some items of (\ref{eqn3.7}):
$$
T_m'(t)e^{\frac{\alpha}{1-\alpha}t}=\tilde{T}_m'(t)-\frac{\alpha}{1-\alpha}\tilde{T}_m(t),
$$
$$
T_m''(t)e^{\frac{\alpha}{1-\alpha}t}=\tilde{T}_m''(t)-\frac{2\alpha}{1-\alpha}\tilde{T}_m'(t)+\left(\frac{\alpha}{1-\alpha}\right)^2\tilde{T}_m(t),
$$
...
\begin{equation} 
T_m^{(n)}(t)e^{\frac{\alpha}{1-\alpha}t}=\sum\limits_{j=0}^n (-1)^{n-j}\frac{n!}{j!(n-j)!}\left(\frac{\alpha}{1-\alpha}\right)^{n-j}\tilde{T}_m^{(j)}(t).
\label{eqn3.8}
\end{equation}
We note that $\tilde{T}_m^{(0)}(t)=\tilde{T}_m(t)$.

Considering (\ref{eqn3.8}), from (\ref{eqn3.7}) we deduce
\begin{equation}
\begin{array}{l}
\sum\limits_{n=0}^k \lambda_n\left\{\sum\limits_{i=0}^n\left(-\frac{\alpha}{1-\alpha}\right)^i \sum\limits_{j=0}^{n-i}\frac{(n-i)!}{j!(n-i-j)!}\left(-\frac{\alpha}{1-\alpha}\right)^{n-i-j}\times\right.\\
\left.\left[\tilde{T}_m^{(j)}(t)-\tilde{T}_m^{(j)}(0)\right]+\left(-\frac{\alpha}{1-\alpha}\right)^{n+1}\int\limits_0^t \tilde{T}_m(s)ds\right\}+(m\pi)^2(1-\alpha)\tilde{T}_m(t)=\\
(1-\alpha)e^{\frac{\alpha}{1-\alpha}t}f_m(t).
\end{array}
\label{eqn3.9}
\end{equation}
Differentiating (\ref{eqn3.9}) once by $t$, we will get $(k+1)$th order differential equation. Using its general solution and applying initial conditions, one can get explicit form of functions $T_m(t)$, consequently, formal solution of the formulated problem, represented by infinite series (\ref{eqn3.4}). Imposing certain conditions to the given functions, we prove uniform convergence of infinite series, which will complete the proof of a unique solvability of the formulated problem.

In the next section we will show complete steps in particular case. We note that even this particular case was not considered before.

\section{Particular case}
In this subsection we will consider particular case $k=2$, in order to show complete steps. In this case, after differentiation of (\ref{eqn3.9}) once by $t$, we will get the following third order ordinary differential equation
\begin{equation}
\tilde{T}_m'''(t)+A_1\tilde{T}_m''(t)+A_2\tilde{T}_m'(t)+A_3\tilde{T}_m=g_m(t)
\label{eqn4.1}
\end{equation}
where $\tilde{T}_m(t)=T_m(t)e^{\frac{\alpha}{1-\alpha}t}$,
$$
A_1=-\frac{3\alpha}{1-\alpha}+\frac{\lambda_1}{\lambda_2} ,\,\,\,A_2=3\left(-\frac{\alpha}{1-\alpha}\right)^2+\frac{\lambda_0-\frac{2\alpha\lambda_1}{1-\alpha}+\left(m\pi\right)^2(1-\alpha)}{\lambda_2},
$$
$$
A_3=\left(-\frac{\alpha}{1-\alpha}\right)^3+\frac{\left(-\frac{\alpha}{1-\alpha}\right)^2\lambda_1-\frac{\alpha}{1-\alpha}\lambda_0}{\lambda_2},
$$
\begin{equation}
g_m(t)=[\alpha{f}_m(t)+(1-\alpha){f_m}'(t)]e^{\frac{\alpha}{1-\alpha}t}.
\label{eqn4.2}
\end{equation}
\subsection{General solution}
Characteristic equation of (\ref{eqn4.1}) is 
$$ 
\mu^3+A_1\mu^2+A_2\mu+A_3=0,
$$  
discriminant of which is 
$$ 
\Delta_m=-4A_1^3A_3+A_1^2A_2^2-4A_2^3+18A_1A_2A_3-27A_3^2.
$$
According to the general theory, form of solutions depend on the sign of the determinant. Below we will give explicit forms of solutions case by case.

Case $\Delta_m>0$. In this case, characteristic equation will have 3 different real roots ($\mu_1,\mu_2,\mu_3$), hence based on general solution we can find explicit form of $T_m(t)$ as
\begin{equation}
\begin{array}{l} 
{T_m(t)=C_{1} e^{(\mu _{1} -\frac{\alpha }{1-\alpha } )t} +C_{2} e^{(\mu _{2} -\frac{\alpha }{1-\alpha } )t} +C_{3} e^{(\mu _{3} -\frac{\alpha }{1-\alpha } )t} +} \\ 
{+\frac{1}{(\mu _{2} \mu _{3}^{2} -\mu _{2}^{2} \mu _{3} -\mu _{1} \mu _{3}^{2} +\mu _{1}^{2} \mu _{3} +\mu _{1} \mu _{2}^{2} -\mu _{1}^{2} \mu _{2} )} \left[e^{(\mu _{1} -\frac{\alpha }{1-\alpha } )t} (\mu _{3} -\mu _{2} )\times\right.}\\
\times\int \left(\alpha g_m(z)+(1-\alpha )g_m'(z)\right)e^{(\frac{\alpha }{1-\alpha } -\mu _{1} )z} dz+ e^{(\mu _{2} -\frac{\alpha }{1-\alpha } )t} (\mu _{1} -\mu _{3} )\times\\
\times\int (\alpha g_m(z)+(1-\alpha )g_m'(z))e^{(\frac{\alpha }{1-\alpha } -\mu _{2} )z} dz+\\
{\left.+e^{\left(\mu _{3} -\frac{\alpha }{1-\alpha } \right)t} (\mu _{2} -\mu _{1} )\int (\alpha g_m(z)+(1-\alpha )g_m'(z))e^{(\frac{\alpha }{1-\alpha } -\mu _{3} )z}   \right]}. 
\end{array}
\label{eqn4.3}
\end{equation}

Case $\Delta_m<0$. In this case, characteristic equation has one real ($\mu_1$) and two complex-conjugate roots ($\mu_2=\mu_{21}\pm i\mu_{22}$). Therefore, $T_m(t)$ will have a form
\begin{equation}
\begin{array}{l} {T_m(t)=C_{1} e^{(\mu _{1} -\frac{\alpha }{1-\alpha } )t} +(C_{2} \cos \mu _{22} t\, +C_{3} \sin \mu _{22} t)e^{(\mu _{21} -\frac{\alpha }{1-\alpha } )t} +} \\
 {+\frac{1}{\mu _{22}^{2} -3\mu _{21}^{2} -2\mu _{1} \mu _{21} +\mu _{1}^{2} } \left[e^{(\mu _{1} -\frac{\alpha }{1-\alpha } )t} \int (\alpha g_m(z)+(1-\alpha )g_m'(z))e^{(\frac{\alpha }{1-\alpha } -\mu _{1} )} dz +\right. } \\
  {+\frac{1}{\mu _{22} } e^{(\mu _{21} -\frac{\alpha }{1-\alpha } )t} \cos \mu _{22} t\int (\mu _{1} \sin \mu _{22} z-\mu _{21} \sin \mu _{22} z-\mu _{22} \cos \mu _{22} z)\times}\\
  (\alpha g_m(z)+(1-\alpha )g_m'(z))e^{(\frac{\alpha }{1-\alpha } -\mu _{21} )z} dz + \\ 
  {+ \frac{1}{\mu _{22} } e^{(\mu _{21} -\frac{\alpha }{1-\alpha } )t} \sin \mu _{22} t  \int (\mu _{21} \cos \mu _{22} z-\mu _{22} \sin \mu _{22} z-\mu _{1} \cos \mu _{22} z)\times}\\
  (\alpha g_m(z)+(1-\alpha )g_m'(z))e^{(\frac{\alpha }{1-\alpha } -\mu _{21} )z} dz ]. 
  \end{array}
\label{eqn4.4}
\end{equation}

Case $\Delta_m=0$. This case will have 2 sub-cases: 

a) 3 real roots, two of which are equal, third one is different ($\mu_1=\mu_2,\,\mu_3$): 
\begin{equation}
\begin{array}{l} 
{T_m(t)=C_{1} e^{(\mu _{1} -\frac{\alpha }{1-\alpha } )t} +C_{2} te^{(\mu _{1} -\frac{\alpha }{1-\alpha } )t} +C_{3} e^{(\mu _{3} -\frac{\alpha }{1-\alpha } )t} +} \\
 {+\frac{1}{(\mu _{1}^{2} -2\mu _{1} \mu _{3} +\mu _{3} )} \left[e^{(\mu _{1} -\frac{\alpha }{1-\alpha } )t} \int (\mu _{3} z-1-\mu _{1} z)(\alpha g_m(z)+(1-\alpha )g_m'(z))e^{(\frac{\alpha }{1-\alpha } -\mu _{1} )z} dz \right. +} \\ 
 {+te^{(\mu _{1} -\frac{\alpha }{1-\alpha } )t} (\mu _{1} -\mu _{3} )\int (\alpha g_m(z)+(1-\alpha )g_m'(z))e^{(\frac{\alpha }{1-\alpha } -\mu _{1} )z} dz +} \\ 
 {+\left. e^{(\mu _{3} -\frac{\alpha }{1-\alpha } )t} \int (\alpha g_m(z)+(1-\alpha )g_m'(z))e^{(\frac{\alpha }{1-\alpha } -\mu _{3} )z} dz \right]}; 
 \end{array}
\label{eqn4.5}
\end{equation}

b) all 3 real roots are the same ($\mu_1=\mu_2=\mu_3$):
\begin{equation}
\begin{array}{l} {T_m(t)=C_{1} e^{(\mu _{1} -\frac{\alpha }{1-\alpha } )t} +C_{2} te^{(\mu _{1} -\frac{\alpha }{1-\alpha } )t} +C_{3} t^{2} e^{(\mu _{1} -\frac{\alpha }{1-\alpha } )t} +}\\
{+\frac{1}{2} e^{(\mu _{1} -\frac{\alpha }{1-\alpha } )t} \int z^{2} (\alpha g_m(z)+(1-\alpha )g_m'(z))e^{(\frac{\alpha }{1-\alpha } -\mu _{1} )z} dz -} \\ {-te^{(\mu _{1} -\frac{\alpha }{1-\alpha } )t} \int z(\alpha g_m(z)+(1-\alpha )g_m'(z))e^{(\frac{\alpha }{1-\alpha } -\mu _{1} )z} dz +} \\ {\frac{1}{2} t^{2} e^{(\mu _{1} -\frac{\alpha }{1-\alpha } )t} \int (\alpha g_m(z)+(1-\alpha )g_m'(z))e^{(\frac{\alpha }{1-\alpha } -\mu _{1} )z} dz }. \end{array}
\label{eqn4.6}
\end{equation}
Here $C_j$ ($j=\overline{1,3}$) are any constants, which will be defined using initial conditions.

\subsection{Convergence part}
We consider the case $\Delta_m>0$ in details. Found solution we satisfy to the initial conditions (\ref{eqn3.5}).Without losing generality, we assume that $\tilde{C}_i=0$ ($i=\overline{0,2}$). Regarding the $C_j$ we will get the following algebraic system of equations
$$
\left\{\begin{array}{l} {C_{1} +C_{2} +C_{3} =-d_{1} } \\ {C_{1} (\mu _{1} -\frac{\alpha }{1-\alpha } )+C_{2} (\mu _{2} -\frac{\alpha }{1-\alpha } )+C_{3} (\mu _{3} -\frac{\alpha }{1-\alpha } )=-d_{2} } \\ {C_{1} (\mu _{1} -\frac{\alpha }{1-\alpha } )^{2} +C_{2} (\mu _{2} -\frac{\alpha }{1-\alpha } )^{2} +C_{3} (\mu _{3} -\frac{\alpha }{1-\alpha } )^{2} =-d_{3} } \end{array}\right.,
$$ 
where
\begin{equation}
\begin{array}{l} 
{d_{1} =\frac{1}{(\mu _{2} \mu _{3}^{2} -\mu _{2}^{2} \mu _{3} -\mu _{1} \mu _{3}^{2} +\mu _{1}^{2} \mu _{3} +\mu _{1} \mu _{2}^{2} -\mu _{1}^{2} \mu _{2} )} \left[(\mu _{3} -\mu _{2} )\int _{}^{t}\alpha g_m(t)+\right.}\\
	\left.+(1-\alpha )g_m'(t)dt\right|_{t=0}  
 +(\mu _{1} -\mu _{3} )\int _{}^{t}\left. \alpha g_m(t)+(1-\alpha )g_m'(t)dt\right|_{t=0} +\\
 \left.+(\mu _{2} -\mu _{1} )\int _{}^{t}\left. \alpha g_m(t)+(1-\alpha )g_m'(t)dt\right|_{t=0}  \right] 
 \end{array}
 \label{eqn4.7}
 \end{equation}
 
  \begin{equation}
 \begin{array}{l} 
 {d_{2} =\frac{1}{(\mu _{2} \mu _{3}^{2} -\mu _{2}^{2} \mu _{3} -\mu _{1} \mu _{3}^{2} +\mu _{1}^{2} \mu _{3} +\mu _{1} \mu _{2}^{2} -\mu _{1}^{2} \mu _{2} )} \left[(\mu _{1} -\frac{\alpha }{1-\alpha } )(\mu _{3} -\mu _{2} )\int _{}^{t}\left. \alpha g_m(t)+(1-\alpha )g_m'(t)dt\right|_{t=0}  \right. +} 
 \\ {+(\mu _{3} -\mu _{2} )(\alpha g_m(0)+(1-\alpha )g_m'(0))+(\mu _{2} -\frac{\alpha }{1-\alpha } )(\mu _{1} -\mu _{3} )\int _{}^{t}\left. \alpha g_m(t)+(1-\alpha )g_m'(t)dt\right|_{t=0}  +} \\ {+(\mu _{1} -\mu _{3} )(\alpha g_m(0)+(1-\alpha )g_m'(0))+(\mu _{3} -\frac{\alpha }{1-\alpha } )(\mu _{2} -\mu _{1} )\int _{}^{t}\left. \alpha g_m(t)+(1-\alpha )g_m'(t)dt\right|_{t=0}  +} \\ {\left. +(\mu _{2} -\mu _{1} )(\alpha g_m(0)+(1-\alpha )g_m'(0))\right]}
  \end{array}
  \label{eqn4.8}
  \end{equation}
  
  \begin{equation}
  \begin{array}{l}
  {d_{3} =\frac{1}{(\mu _{2} \mu _{3}^{2} -\mu _{2}^{2} \mu _{3} -\mu _{1} \mu _{3}^{2} +\mu _{1}^{2} \mu _{3} +\mu _{1} \mu _{2}^{2} -\mu _{1}^{2} \mu _{2} )} \left[(\mu _{1} -\frac{\alpha }{1-\alpha } )^{2} (\mu _{3} -\mu _{2} )\int _{}^{t}\left. \alpha g_m(t)+(1-\alpha )g_m'(t)dt\right|_{t=0}  \right. +} \\ 
  {+(\mu _{1} -\frac{\alpha }{1-\alpha } )(\mu _{3} -\mu _{2} )(\alpha g_m(0)+(1-\alpha )g_m'(0))+(\mu _{3} -\mu _{2} )(\alpha g_m'(0)+(1-\alpha )g_m''(0))+} \\ 
  {(\mu _{2} -\frac{\alpha }{1-\alpha } )^{2} (\mu _{1} -\mu _{3} )\int _{}^{t}\left. \alpha g_m(t)+(1-\alpha )g_m'(t)dt\right|_{t=0}  +(\mu _{2} -\frac{\alpha }{1-\alpha } )(\mu _{1} -\mu _{3} )(\alpha g_m(0)+(1-\alpha )g_m'(0))+} \\ 
  {+(\mu _{1} -\mu _{3} )(\alpha g_m'(0)+(1-\alpha )g_m''(0))+(\mu _{3} -\frac{\alpha }{1-\alpha } )^{2} (\mu _{2} -\mu _{1} )\int _{}^{t}\left. \alpha g_m(t)+(1-\alpha )g_m'(t)dt\right|_{t=0}  +} \\ 
  {\left. +(\mu _{3} -\frac{\alpha }{1-\alpha } )(\mu _{2} -\mu _{1} )(\alpha g_m(0)+(1-\alpha )g_m'(0))+(\mu _{2} -\mu _{1} )(\alpha g_m'(0)+(1-\alpha )g_m''(0))\right].} 
  \end{array}
  \label{eqn4.9}
  \end{equation}

Solving this system, we get
$$
\begin{array}{l}
{C_{1} =-d_{1}-\frac{-d_{2} (\mu _{1} -\frac{\alpha }{1-\alpha } )+d_{3}}{((\mu _{2} -\frac{\alpha }{1-\alpha } )(\mu _{1} -\frac{\alpha }{1-\alpha } )-(\mu _{2} -\frac{\alpha }{1-\alpha } )^{2} )} -} 
\\ {-\frac{-d_{1}(\mu _{1} -\frac{\alpha }{1-\alpha } )(\mu _{2} -\frac{\alpha }{1-\alpha } )-d_{2} (\mu _{2} -\frac{\alpha }{1-\alpha } )+d_{2} (\mu _{1} -\frac{\alpha }{1-\alpha } )-d_{3} }{\left[\mu _{1} \mu _{2} -\mu _{2} \mu _{3} -\mu _{3} \mu _{1} -\mu _{3}^{2} -2(\frac{\alpha }{1-\alpha }) ^{2} \right]\left[((\mu _{2} -\frac{\alpha }{1-\alpha } )(\mu _{1} -\frac{\alpha }{1-\alpha } )-(\mu _{2} -\frac{\alpha }{1-\alpha } )^{2} )\right]} -} 
\\ {-\frac{-d_{1} (\mu _{1} -\frac{\alpha }{1-\alpha } )(\mu _{2} -\frac{\alpha }{1-\alpha } )+d_{2} (\mu _{2} -\frac{\alpha }{1-\alpha } )+d_{2} (\mu _{1} -\frac{\alpha }{1-\alpha } )-d_{3} }{\mu _{1} \mu _{2} -\mu _{2} \mu _{3} -\mu _{3} \mu _{1} -\mu _{3}^{2} -2(\frac{\alpha }{1-\alpha } )^{2} } } \\ {C_{2} =\frac{-d_{2} (\mu _{1} -\frac{\alpha }{1-\alpha } )+d_{3} }{((\mu _{2} -\frac{\alpha }{1-\alpha } )(\mu _{1} -\frac{\alpha }{1-\alpha } )-(\mu _{2} -\frac{\alpha }{1-\alpha } )^{2} )} -} \\ {-\frac{(-d_{1} )(\mu _{1} -\frac{\alpha }{1-\alpha } )(\mu _{2} -\frac{\alpha }{1-\alpha } )+d_{2} (\mu _{2} -\frac{\alpha }{1-\alpha } )+d_{2} (\mu _{1} -\frac{\alpha }{1-\alpha } )-d_{3} }{\left[\mu _{1} \mu _{2} -\mu _{2} \mu _{3} -\mu _{3} \mu _{1} -\mu _{3}^{2} -2(\frac{\alpha }{1-\alpha } )^{2} \right]\left[((\mu _{2} -\frac{\alpha }{1-\alpha } )(\mu _{1} -\frac{\alpha }{1-\alpha } )-(\mu _{2} -\frac{\alpha }{1-\alpha } )^{2} )\right]} } \\ {C_{3} =\frac{(-d_{1} )(\mu _{1} -\frac{\alpha }{1-\alpha } )(\mu _{2} -\frac{\alpha }{1-\alpha } )+d_{2} (\mu _{2} -\frac{\alpha }{1-\alpha } )+d_{2} (\mu _{1} -\frac{\alpha }{1-\alpha } )-d_{3} }{\mu _{1} \mu _{2} -\mu _{2} \mu _{3} -\mu _{3} \mu _{1} -\mu _{3}^{2} -2(\frac{\alpha }{1-\alpha } )^{2} } } \end{array}
$$
In general, we can write
$$
\left|C_j\right|\leq M_1|d_1|+M_2|d_2|+M_3|d_3|.
$$
Hence, we need the following estimations in order to provide convergence of used series:
\begin{equation}
\begin{array}{l} {\left|d_{1} \right|\le \frac{1}{(m\pi )^{4} } \left|M_{4} \right. \int (\alpha  f_{4,0} (z)+(1-\alpha )f_{4,1} (z))e^{(\frac{\alpha }{1-\alpha } -\mu _{1} )z} \left. dz\right|_{t=0} +} \\ {+M_{5} \int (\alpha  f_{4,0} (z)+(1-\alpha )f_{4,1} (z))e^{(\frac{\alpha }{1-\alpha } -\mu _{2} )z} \left. dz\right|+} \\ {+M_{6} \int (\alpha  f_{4,0} (z)+(1-\alpha )f_{4,1} (z))e^{(\frac{\alpha }{1-\alpha } -\mu _{3} )z} \left. dz\right|\left. _{t=0} \right|\le \frac{1}{(m\pi )^{4} } M_{7} } \end{array} 
\label{eqn4.10}
\end{equation}
\begin{equation}
\begin{array}{l} {\left|d_{2}\right|\le \frac{1}{(m\pi )^{4} } \left|M_{8} \right. e^{(\mu _{1} -\frac{\alpha }{1-\alpha } )t} \int (\alpha  f_{4,0} (z)+(1-\alpha )f_{4,1} (z))e^{(\frac{\alpha }{1-\alpha } -\mu _{1} )z} \left. dz\right|_{t} +} \\ 
{+M_{9} e^{(\mu _{2} -\frac{\alpha }{1-\alpha } )t} \int (\alpha  f_{4,0} (z)+(1-\alpha )f_{4,1} (z))e^{(\frac{\alpha }{1-\alpha } -\mu _{2} )z} \left. dz\right|_{t} +} \\ 
+M_{10} e^{(\mu _{3} -\frac{\alpha }{1-\alpha } )t} \int (\alpha  f_{4,0} (z)+(1-\alpha )f_{4,1} (z))e^{(\frac{\alpha }{1-\alpha } -\mu _{3} )z} \left. dz\right|_{t} +\\
\left.+M_{11} (\alpha f_{4,0} (z)+(1-\alpha )f_{4,1} (z))\right|\le \frac{1}{(m\pi )^{4} } M_{12}  \end{array}
\label{eqn4.11} 
\end{equation}

\begin{equation}
\begin{array}{l} {\left|d_{3}\right|\le \frac{1}{(m\pi )^{4} } \left|M_{13} \right. e^{(\mu _{1} -\frac{\alpha }{1-\alpha } )t} \int (\alpha  f_{4,0} (z)+(1-\alpha )f_{4,1} (z))e^{(\frac{\alpha }{1-\alpha } -\mu _{1} )z} \left. dz\right|_{t} +} \\ {+M_{14} e^{(\mu _{2} -\frac{\alpha }{1-\alpha } )t} \int (\alpha  f_{4,0} (z)+(1-\alpha )f_{4,1} (z))e^{(\frac{\alpha }{1-\alpha } -\mu _{2} )z} \left. dz\right|_{t} +} \\ {+M_{15} e^{(\mu _{3} -\frac{\alpha }{1-\alpha } )t} \int (\alpha  f_{4,0} (z)+(1-\alpha )f_{4,1} (z))e^{(\frac{\alpha }{1-\alpha } -\mu _{3} )z} \left. dz\right|_{t} +} \\ {\left. +M_{16} (\alpha f_{4,0} (z)+(1-\alpha )f_{4,1} (z))+M_{17} (\alpha f_{4,1} (z)+(1-\alpha )f_{4,2} (z))\right|\le \frac{1}{(m\pi )^{4} } M_{18} } \end{array} 
\label{eqn4.12}
\end{equation}
Here $M_i\,(i=\overline{1,18})$ are any positive constants,
\begin{eqnarray}
f_m(t)=\int _{0}^{1}f(t,x)\sin m\pi xdx=\frac{1}{(m\pi )^{4} }  f_{4,0} (t),\\
f_m'(t)=\frac{1}{(m\pi )^{4} } \int _{0}^{1}\frac{\partial }{\partial {t} } \left(\frac{\partial ^{4} }{\partial {x}^{4} } f(t,x)\right)\sin m\pi xdx=\frac{1}{(m\pi )^{4} }  f_{4,1} (t),\\ 
f_m''(t)=\frac{1}{(m\pi )^{4} } \int _{0}^{1}\frac{\partial }{\partial^2 {t} } \left(\frac{\partial ^{4} }{\partial {x}^{4} } f(t,x)\right)\sin m\pi xdx=\frac{1}{(m\pi )^{4} }  f_{4,2} (t),
\label{eqn4.13}
\end{eqnarray} 
\begin{eqnarray}
f_{4,0} (t)=\int _{0}^{1}\frac{\partial ^{4} f(t,x)}{\partial {x}^{4} }  \sin m\pi xdx,
\\
f_{4,1} (t)=\, \int _{0}^{1}\frac{\partial ^{5} f(t,x)}{\partial {t} \partial {x}^{4} }  \sin m\pi xdx,
\\
f_{4,2} (t)=\, \int _{0}^{1}\frac{\partial ^{6} f(t,x)}{\partial {t}^{2} \partial {x}^{4} }  \sin m\pi xdx.
\label{eqn4.14}
\end{eqnarray}
We note that for above-given estimations,we need to impose  the following conditons to the given function $f(t,x)$:
\begin{equation}
\begin{array}{l}
\displaystyle{\frac{\partial{f}}{\partial{t}}|_{t=0}=0, \,\,\,\,\, \frac{\partial^2{f}}{\partial{t}^2}|_{t=0}=0, \,\,\,\, \frac{\partial^3{f}}{\partial{t}^3}|_{t=0}=0,}\\
\displaystyle{f(t,1)=f(t,0)=0,\,\,\frac{\partial ^{2} f(t,1)}{\partial{x}^{2} } =\frac{\partial ^{2} f(t,0)}{\partial {x}^{2} } =0.}
\end{array}\label{eqn4.15}
\end{equation}

Based on estimations (\ref{eqn4.7})-(\ref{eqn4.9}), we obtain
\begin{equation}
\left|C_j\right|\leq \frac{M_{19}}{(m\pi)^4}
\label{eqn4.16}
\end{equation}
and considering (\ref{eqn4.3}), finally we get
\begin{equation}
\left|T_m\right|\leq \frac{M_{20}}{(m\pi)^4}.
\label{eqn4.17}
\end{equation}

Taking (\ref{eqn3.4}) into account, on can easily deduce that
\begin{equation}
\left|u(t,x)\right|\leq \frac{M_{21}}{(m\pi)^4},\,\,\,\left|u_{xx}(t,x)\right|\leq \frac{M_{22}}{(m\pi)^2}.
\label{eqn4.18}
\end{equation}

Another required estimation 
\begin{equation}
\left|{}_{CF}D_{0t}^\alpha u(t,x)\right|\leq \frac{M_{23}}{(m\pi)^4}
\label{eqn4.19}
\end{equation}
can be deduced easily, as well.

{\bf Theorem.}
	If $f(t,x)\in C^2(\overline{\Omega})$, $\frac{\partial^3 f(t,x)}{\partial t^3}\in C(\Omega\cup\{t=0\})$ together with (\ref{eqn4.15}), then problem (\ref{eqn3.1})-(\ref{eqn3.3}), when $k=2$ has a unique solution represented by (\ref{eqn3.4}), where $T_m(t)$ are defined by (\ref{eqn4.3})-(\ref{eqn4.6}) depending on the sign of $\Delta_m$.
	
	{\bf Remark 1.} Similar result can be obtained for general case, as well.
	
	{\bf Remark 2.} We note that used algorithm allows us to investigate fractional spectral problems such
	\begin{equation*}
	\left\{
	\begin{array}{l}
	{}_{CF}D_{0t}^{\alpha+1}T(t)+\mu T(t)=0,\\
	T(0)=0,\,\,T(1)=0,\\
	\end{array}
	\right.
	\end{equation*}
	reducing it to the second order usual spectral problem.



\begin{thebibliography}{99}
	
\bibitem{CF1} M. Caputo and M. Fabrizio, A new definition of fractional derivative without singular kernel, Progr. Frac. Differ. Appl. 1(2), 1–13 (2015). 

\bibitem{CF2} M. Caputo and M. Fabrizio, Applications of New Time and Spatial Fractional Derivatives with Exponential Kernels, Progr. Frac. Differ. Appl. 2(1), 1–11 (2016).
 
\bibitem{NL} J. Losada, J.J. Nieto, Properties of a New Fractional Derivative without Singular Kernel, Progr. Frac. Differ. Appl. 1(2), 87–92 (2015). 

\bibitem{AA1} A. Atangana, B.S.T. Alkahtani, New model of groundwater flowing within a confine aquifer: application of Caputo-Fabrizio derivative, Arab J. Geosci. 9: 8 (2016) DOI 10.1007/s12517-015-2060-8

\bibitem{AB1}  A. Atangana, D. Baleanu, Caputo-Fabrizio Derivative Applied to Groundwater Flow within Confined Aquifer,  J. Eng. Mech., D4016005, 1-5 (2016) DOI: 10.1061/(ASCE)EM.1943-7889.0001091 

\bibitem{AA2} A. Atangana, B.S.T. Alkahtani, Extension of the resistance, inductance, capacitance electrical circuit to fractional derivative without singular kernel, Advances in Mechanical Engineering, 7(6), 1–6 (2015) DOI: 10.1177/1687814015591937 

\bibitem{AA3} B.S.T. Alkahtani, A. Atangana, Controlling the wave movement on the surface of shallow water with the Caputo–Fabrizio derivative with fractional order, Chaos, Solitons and Fractals, 89, 539-546 (2016) DOI: 10.1016/j.chaos.2016.03.012  

\bibitem{A1} A. Atangana, On the new fractional derivative and application to nonlinear Fisher’s reaction–diffusion equation, Applied Mathematics and Computaton, 273, 948-956 (2016) DOI: 10.1016/j.amc.2015.10.021

\bibitem{GA1} J.F. G\'{o}mez-Aguilar et al, Modeling of a Mass-Spring-Damper System by Fractional Derivatives with and without a Singular Kernel, Entropy, 17, 6289-6303 (2015) DOI:10.3390/e17096289

\bibitem{GA2} V.F. Morales-Delgado et al, On the solutions of fractional order of evolution equations, Eur. Phys. J. Plus, 132:47 (2017) DOI: 10.1140/epjp/i2017-11341-0

\bibitem{GY} F. Gao et al, Fractional Maxwell fluid with fractional derivative without singular kernel, Thermal Science, 20(S3), S871-S877 (2016)

\bibitem{ASh} Farhad Ali et al, Application of Caputo-Fabrizio derivatives to MHD free convection flow of generalized Walters’-B fluid model, Eur. Phys. J. Plus, 131: 377 (2016) DOI 10.1140/epjp/i2016-16377-x

\bibitem{AN} A. Atangana, J.J. Nieto, Numerical solution for the model of RLC circuit via the fractional derivative without singular kernel, Advances in Mechanical Engineering, 7(10), 1-7 (2015) DOI: 	10.1177/1687814015613758

\bibitem{GA3} J.F. G\'{o}mez-Aguilar et al, Modeling diffusive transport with a fractional derivative without singular kernel, Physica A, 447, 467-481 (2016) DOI: 10.1016/j.physa.2015.12.066  

\bibitem{H}  J. Hristov, Transient heat diffusion with a non-singular fading memory, Thermal Science, 20(2), 757-762 (2016) 

\bibitem{K1} N. Al-Salti, E. Karimov and S. Kerbal, Boundary-value problem for fractional heat equation involving Caputo-Fabrizio derivative, New Trend in Mathematical Sciences, 4, 79-89 (2016)

\bibitem{K2} N. Al-Salti, E. Karimov and K. Sadarangani , On a Differential Equation with Caputo-Fabrizio Fractional Derivative of Order $1<\beta \leq 2$ and Application to Mass-Spring-Damper System, Progr. Fract. Differ. Appl., 2(4), 257-263 (2016)

\end{thebibliography}
\end{document}